\documentclass[11pt]{amsart}
\usepackage{graphics}
\usepackage{amsmath,amsthm,amscd,amssymb}
\usepackage{epsfig}
\usepackage[all,cmtip]{xy}
\usepackage{tikz}
\usetikzlibrary{arrows}
\usetikzlibrary{calc}

\renewcommand{\P}{\mathcal{P}}

\newcommand{\T}{T}

\renewcommand{\d}{\partial}

\newcommand{\D}{\Delta}
\renewcommand{\L}{\Lambda}

\renewcommand{\P}{\mathcal{ P}}
\newcommand{\oM}{\overline{\mathcal{M}}}
\newcommand{\M}{\mathcal{M}}

\newcommand{\rad}{\operatorname{rad}}
\newcommand{\soc}{\operatorname{soc}}

\newtheorem{Theorem}{Theorem}[section]
\newtheorem{Lemma}[Theorem]{Lemma}
\newtheorem{Corollary}[Theorem]{Corollary}
\newtheorem{Definition}[Theorem]{Definition}
\theoremstyle{definition}
\newtheorem*{remark}{Remark}

\begin{document}

\topmargin .7cm \oddsidemargin 1.5cm \evensidemargin 1.5cm

\title[Trivial Extensions of Gentle Algebras]{Trivial Extensions of Gentle Algebras and Brauer Graph Algebras}
\thanks{This work was supported through the Engineering and Physical Sciences Research Council, grant number EP/K026364/1, UK}

\author{Sibylle Schroll}
\address{Sibylle Schroll\\
Department of Mathematics \\
University of Leicester \\
University Road \\
Leicester LE1 7RH \\
United Kingdom}
\email{schroll@le.ac.uk}


\begin{abstract}
We show that two well-studied classes of tame algebras coincide: namely, the class of symmetric special biserial algebras coincides with the class of Brauer graph algebras. We then explore the connection between gentle algebras and symmetric special biserial algebras by explicitly determining 
the trivial extension of a gentle algebra by its minimal injective co-generator. This is a symmetric special biserial algebra and hence a Brauer graph algebra of which we explicitly give the Brauer graph. 
We further show that a Brauer graph algebra gives rise,  via admissible cuts, to many  gentle algebras and that the trivial extension of a gentle algebra obtained via an admissible cut is the original Brauer graph algebra. 

As a consequence we prove that the trivial extension of a Jacobian algebra of an ideal triangulation of a Riemann surface with marked points in the boundary is isomorphic to the Brauer graph algebra with Brauer graph given by the 
arcs of the triangulation. 
\end{abstract}


\maketitle



\section{Introduction}

\parskip7pt
\parindent0pt

One of the cornerstones of the representation theory of finite dimensional algebras is Drozd's classification of algebras in terms of their representation type, which is either 
 finite, tame or wild. Algebras of tame representation type are of great interest since they have infinitely many isomorphism classes of indecomposable modules, yet they do usually still exhibit discernible patterns in their representation theory.

One of the most important classes of tame algebras are special biserial algebras. 
This class comprises many well-studied classes of algebras such as gentle algebras, string algebras and Brauer graph algebras. 



In this paper we start by proving that the classes of Brauer graph algebras and that of symmetric special biserial
 algebras coincide. We then determine that the trivial extension of a gentle algebra by its minimal injective co-generator is a Brauer graph algebra, 
 for which we explicitly construct the Brauer graph. 
%
%

While trivial extensions are a way to obtain new algebras from existing ones by increasing the dimension of the underlying vector space, admissible cuts, 
as introduced in \cite{F}, are a way of obtaining new algebras from existing ones by decreasing the dimension of the underlying vector
space. 
We show that in the setting of symmetric special biserial and gentle algebras, taking admissible cuts and taking trivial extensions are 
inverse operations in the following sense. We show that an admissible cut of a symmetric special biserial algebra (with simple cycles) or equivalently of a Brauer graph algebra (without multiplicities) is a gentle algebra and that the original 
algebra is recovered from this gentle algebra via a trivial extension. 

This is of particular interest since it can be directly applied to Jacobian 
algebras of triangulations of Riemann surfaces with marked points in the boundary to show that the trivial extension of  a 
Jacobian algebra is the Brauer 
graph algebra with Brauer graph given by the arcs of the triangulation. 


Historically, the study of Brauer graph algebras was initiated by Richard Brauer in the context of modular representation 
theory of finite groups. They play a pivotal role as they appear in the form of Brauer tree algebras as blocks of finite groups with cyclic defect group. 
However, since their inception, Brauer graph algebras have also been extensively studied outside the context of the representation theory of finite groups, 
see for example  \cite{GSST, K, Ric, Ro}, for a small selection of this literature. In parallel and largely independently, 
special biserial algebras and in particular, 
symmetric special biserial algebras have been extensively studied. For example, their Auslander-Reiten structure has been determined in \cite{ES},   
 and, for an example of more recent papers on the cohomological structure of these algebras see \cite{E, ESch, ST}.

\sloppy Roggenkamp has shown in \cite{Ro}  that  if the quiver of a symmetric special biserial algebra $\Lambda$ has no two arrows in the same direction
then $\Lambda$ is a Brauer graph algebra. We show that this result holds in general.  That is 
any symmetric special biserial algebra is a Brauer graph algebra. We do this by associating a graph with local structure to a
symmetric special biserial algebra and we  show that if this graph is considered as a Brauer graph  then the 
corresponding Brauer graph algebra is isomorphic to the 
symmetric special biserial algebra. 
More precisely,  we show the following.


\begin{Theorem}\label{SymmetricSpecialBiserialIsBGA} 
Let $k$ be an algebraically closed field and let $\Lambda = kQ/I$ be a finite dimensional symmetric special biserial algebra and let $G_\Lambda$ be its graph with local structure. 
Let $B$ be the  Brauer graph algebra associated to $G_\Lambda$. Then $B$ is isomorphic to $\Lambda$. 
\end{Theorem}

Given a gentle algebra $A$, the trivial extension $T(A) = A \ltimes D(A)$ of $A$ by its minimal co-generator $D(A)$ is a symmetric special biserial algebra \cite{PS, Ringel, Schroer}. Thus it is a Brauer 
graph algebra. In order to identify its Brauer graph we construct a graph $\Gamma_A$, equipped with a cyclic ordering of the edges around each vertex induced by the maximal paths 
in $A$ and we show that $\Gamma_A$ is the Brauer graph of $T(A)$.  

\begin{Theorem} \label{TrivialExtensionofGentleAlgebra}
Let $k$ be an algebraically closed field and let $A= kQ/I$ be a finite dimensional gentle algebra and let $\Gamma_A$ be its graph. Let $B$ be the Brauer graph algebra defined on $\Gamma_A$ with cyclic ordering induced by the maximal paths in $A$. Then the trivial extension $T(A)$ of $A$ is isomorphic to $B$ and  $G_{T(A)} = \Gamma_A$ where $G_{T(A)}$ is the graph of $T(A)$ as symmetric special biserial algebra. 
\end{Theorem}

We show that every Brauer graph algebra with multiplicity identically one (or alternatively every symmetric special biserial algebra in which all cycles in the relations are no power
of a proper sub-cycle) is the trivial extension of a gentle algebra. For this we define the notion of an admissible cut of a Brauer graph algebra. Admissible cuts of this form were first
considered in \cite{F, FP2}.

\begin{Theorem}\label{ArrowDeletionTheorem}
Let $k$ be an algebraically closed field and let $\Lambda = kQ_\Lambda / I_\Lambda$ be a  Brauer graph algebra with multiplicity one at all vertices in the associated Brauer graph.
 Let $A$ be an admissible cut of $\Lambda$. 
Then $A$ is gentle and $T(A)$ and $\Lambda$ are isomorphic. 
\end{Theorem}

In fact, there are many different gentle algebras that can be obtained from a given Brauer graph algebra resulting from different admissible cuts. In general, while having the same number of isomorphism classes of simple modules these algebras are non-isomorphic and usually they are not even derived equivalent. However, as shown in 
Theorem~\ref{ArrowDeletionTheorem}, their trivial extensions are all isomorphic.  

Furthermore, we have the following consequence of Theorem~\ref{ArrowDeletionTheorem}.

\begin{Corollary}\label{BGAisTrivialExtension}
Every Brauer graph algebra with multiplicity function identically one is the trivial extension of a gentle algebra.
\end{Corollary}

Connections between  Brauer graph algebras and Jacobian algebras of triangulations of marked Riemann surfaces  have recently been 
established by several authors. In particular, the connection of mutation,  flip of diagonals in triangulations and derived equivalences 
have been studied  in \cite{Lad1, Lad2}, \cite{Ai} and \cite{MS}. In particular, a flip of a diagonal always gives rise to a derived equivalence 
of the corresponding Brauer graph algebras \cite{MS}, however, the corresponding Jacobian algebras are not necessarily derived equivalent \cite{Lad1}.

Here we give another connection of Brauer graph algebras and Jacobian algebras of marked surfaces. Namely, 
let $(S,M)$ be a bordered Riemann surface and $M$ a set of marked points in the boundary of $S$. Given an ideal triangulation $\T$ of $(S,M)$, 
let $(Q,I)$ be the associated bound quiver as
defined in \cite{CCS, DWZ} and let $A = kQ/I$ be the associated finite dimensional gentle algebra \cite{ABCP}. 
The orientation of $S$ induces a cyclic ordering of the arcs of $\T$ around each marked point and thus $\T$ can be viewed as a Brauer graph. 

\begin{Corollary}\label{Jacobian}
Let $A$ be a gentle algebra over an algebraically closed field arising from an ideal triangulation $\T$ of a Riemann surface with a set of marked points in the boundary.  Then the trivial extension $T(A)$ of $A$ is isomorphic to 
the Brauer graph algebra with Brauer graph $\T$. 
\end{Corollary}

Given a triangulation $\T$ of  a   Riemann surface with marked points in the boundary $(S,M)$, in \cite{DS}  
surface algebras were defined by cutting $\T$ at internal triangles. This gives
rise to a partial triangulation $P$ of $(S,M)$. A surface algebra $A$
is  a finite dimensional gentle algebra and it can be constructed by associating a bound quiver $(Q,I)$ to $P$
such 
that $A = kQ/I$, see \cite{DS}.  We show that  $P$ together with the cyclic ordering of arcs around each marked point
 induced by the orientation of the surface
  is the Brauer graph associated to the 
trivial extension of $A$. 

\begin{Corollary}\label{SurfaceAlgebra}
Let $k$ be an algebraically closed field and let the $k$-algebra $A$ be a surface algebra of a partial triangulation $P$ of a  Riemann surface with a set of marked points in the boundary. 
Then the trivial extension $T(A)$ of $A$ is isomorphic to 
the Brauer graph algebra with Brauer graph $P$. 
\end{Corollary}

\textbf{Acknowledgements.} The author would like to thank Robert Marsh for helpful conversations and the careful reading of the first draft 
of the present paper, \O yvind Solberg 
for the provision of some examples of trivial extensions, and Ilke Canakci for  an initiation to TikZ pictures. The author also would like to thank the referee for their helpful comments.

\section{Gentle algebras, symmetric special biserial algebras and Brauer graph algebras} 

Let $k$ be an algebraically closed field and let $Q$ be a finite connected quiver. 
Let $I$ be an admissible ideal in the path algebra $kQ$ such that $kQ/I$ is a finite dimensional algebra.
A path in $Q$ is in the bound quiver $(Q,I)$ if it avoids the relations in $I$.

 For a finite dimensional $k$-algebra $A$, let $D= Hom_k(-, k)$ denote the standard duality of   
 the module category $A-mod$ of finitely generated $A$-modules. 
The algebra $A$ is symmetric if it is isomorphic to $D(A)$ as an $A$-$A-$bimodule.
 Let $A^e \simeq A \otimes_k A^{op}$ be the enveloping algebra of $A$. 
 The socle $\soc(M)$ of a right $A$-module $M$ is the largest semisimple submodule of $M$ and the radical of $M$ is defined by $\rad(M) = \rad(A) M $ where $\rad(A)$ is the Jacobson radical of $A$. Unless otherwise stated all algebras considered are indecomposable and all modules are right modules.

\subsection{Special biserial and gentle algebras}

We say that a finite dimensional algebra $A$ is \textit{special biserial} if it is 
Morita equivalent to an algebra of the form $kQ/I$ where

(S1) Each vertex of $Q$ is the starting point of at most two arrows and is the end point of at most two arrows. 

(S2) For each arrow $\alpha$ in $Q$ there is at most one arrow $\beta$ in $Q$ such that $\alpha\beta$ is not in $I$ and there is at most one arrow $\gamma$ such that $\gamma\alpha$ is not in 
$I$. 

We say that an algebra $A$ is \textit{gentle} if it is Morita equivalent to an algebra $kQ/I$ satisfying (S1), (S2) and 

(S3) $I$ is generated by paths of length 2.

(S4) For each arrow $\alpha$ in $Q$ there is at most one arrow $\delta$ in $Q$ such that $\alpha \delta$ is in $I$ and there is at most one   arrow $\varepsilon$ in $Q$ such that $\varepsilon \alpha$ is in $I$. 

An arrow $\alpha$ in $Q$ starts at the vertex $s(\alpha)$ and ends at the vertex $t(\alpha)$. If $p = \alpha_1 \alpha_2 \ldots \alpha_n$, for arrows $\alpha_i$, $1 \leq i \leq n$, is a path in $Q$ then 
$s(p) = s(\alpha_1)$ and $t(p) = t(\alpha_n)$. 


 


Let $\P$ be the set of paths in $(Q,I)$. Let $\M$ be the set of maximal elements in $\P$, that is all paths  $p$ in $\P$ 
such that for all arrows $\alpha$ in $Q$, $\alpha p \notin  \P$ and 
$ p  \alpha \notin  \P$.
In general, for any finite dimensional algebra every path is a subpath of some maximal path and for a gentle algebra that maximal path is unique.  
Furthermore, a non-trivial path $p$ in $Q$ is in $\P$ 
if and only if it is a subpath of a (unique if $A$ gentle) maximal path $m \in \M$ where $m = q p q'$ with $q, q' $ paths in $Q$. 

Suppose $kQ/I$ is a gentle algebra. Then every arrow is contained in a unique maximal path, and there are at most two maximal paths starting at any given vertex, 
and at most two maximal paths ending at any given vertex. It follows that two distinct maximal paths cannot have a common arrow.  Hence  
maximal paths only intersect at a vertex of $Q$. 


\begin{Lemma}\label{MaximalPathsBasis}
Let $ A = kQ/I$ be a gentle algebra. Then the maximal paths in $(Q,I)$ form a basis of $\rm{soc}_{A^e}  A$. 
\end{Lemma}

{\it Proof:} We have ${\rm soc}_{A^e}  A = \bigoplus_{i,j  \;\; vertices \;\; in \;\; Q}  e_i ({\rm soc}_{A^e} A) e_j$ as $k$-vector spaces where 
$e_i$ and  $e_j$  are the trivial paths at vertices $i$ and $j$ respectively.  
We will show that the maximal paths from $i$ to $j$ form a basis of $ e_i ({\rm soc}_{A^e} A) e_j$. Recall that a path $p$ in $(Q,I)$ is maximal
if for all arrows $\alpha$ in $Q$, $p\alpha, \alpha p \in I$ or equivalently if $p\alpha = 0 = \alpha p $ in $A$. Set $R={\rm soc}_{A^e} A$. 
Since $A$ is gentle there are at most two non-zero paths from $i$ to $j$, denote them by $p= \alpha_1 \alpha_2 \ldots \alpha_m$ and $q= \beta_1 \beta_2 \ldots \beta_n$.  Let $\alpha, \beta, \gamma, \delta $ be arrows in $Q$ such that $t(\alpha) = t(\beta) = i $ and $s(\gamma) = s(\delta) = j$ as in the diagram below. 

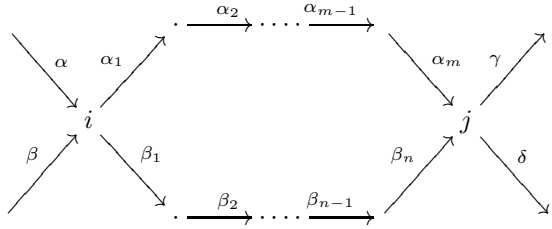
\nopagebreak

{\small

\begin{displaymath}
    \xymatrix{
         \ar[dr]^\alpha  & &\cdot \ar[r]^{\alpha_2} & \cdots \ar[r]^{\alpha_{m-1}} \cdot & \ar[dr]^{\alpha_m}  && \\
        & i  \ar[ur]^{\alpha_1} \ar[dr]^{\beta_1} &&&& j \ar[ur]^{\gamma} \ar[dr]^{\delta}\\
  \ar[ur]^\beta   &&    \cdot \ar[r]^{\beta_2} & \cdots \ar[r]^{\beta_{n-1}} \cdot & \ar[ur]^{\beta_n} &&  }
\end{displaymath}
}

A case by case analysis will give a basis of $e_iRe_j$ in terms of $p$ and $q$. 
Suppose first that none of the arrows $\alpha, \beta, \gamma,  \delta$ exist in $Q$, that is the only arrows starting at $i$ are $\alpha_1$ and $\beta_1$ and there are no arrows ending at $i$ and the only arrows ending at $j$ are $\alpha_m$ and $\beta_n$ and there are no arrows starting at $j$. Then $p$ and $q$ are maximal and $\{p,q\}$ is a basis of $e_iRe_j$. 

Now suppose that  $\beta,  \gamma, \delta$ do not exist in $Q$  but that  $\alpha$ is an arrow in $Q$.   Then either $\alpha p \notin I$ or $\alpha q \notin I $. Suppose that $\alpha p \notin I$. 
Then $p$ is not maximal and $q$ is maximal and $\{q\}$ is a basis of $e_iRe_j$.  
Similarly, if $ \alpha, \beta, \gamma$ do not exist in $Q$ but $\delta$ is an arrow in $Q$.

Suppose now that $\beta$ and $ \delta $ do not exist in $Q$ but that $\alpha$ and $\gamma$ are arrows in $Q$. Then there are two cases: (i) if $\alpha \alpha_1 \in I$ and $\alpha_m \gamma \in I$ then $\{p\}$ is a basis of $e_iRe_j$, (ii) if $\alpha \alpha_1 \in I$ and $\beta_n \gamma \in I$ then $p$ and $q$ are not maximal and no linear combination of $p$ and $q$ is such that it is zero by left and right multiplication by all arrows in $Q$. Hence $e_iRe_j =\{ 0 \}$.

If $\alpha$ and $\beta$ are arrows in $Q$ and $\gamma$ but  $ \delta$ do not exist in $Q$  then either (i) $\alpha \alpha_1 \in I $ and $\beta \beta_1 \in I$ or (ii) $\alpha \beta_1 \in I$ and $\beta \alpha_1 \in I$. Assume without loss of generality that (i) holds. Then $p$ and $q$ are not maximal and as above $e_iRe_j =\{ 0 \}$.

Suppose that $\alpha$, $\beta$, and $\gamma $ are arrows in $Q$ and that $\delta $ does not exist. Then as above we can assume without loss of generality that $\alpha_m \gamma \in I$ and hence that $p$ is maximal but $q$ is not and $\{p\}$ is a basis of $e_i R e_j$. 

If all arrows $\alpha, \beta, \gamma$ and $\delta$ exist in $Q$ then both $p$ and $q$ cannot be maximal and $e_i Re_j = \{0\}$. 

All remaining cases are covered by similar arguments. \hfill $\Box$



\subsection{Brauer graph algebras}\label{BGA}

 
In this section we define symmetric Brauer graph algebras. 

We call a finite connected graph $\Gamma$ (possibly containing loops and multiple edges) a \textit{Brauer graph} if $\Gamma$ is equipped with a cyclic ordering of the edges around each vertex and if for every vertex $\nu$ in  $\Gamma$ there is 
an associated strictly positive integer $e(\nu)$ called the \emph{multiplicity of $\nu$}.

For an edge $E$ in $\Gamma$  with vertices $\nu$ and $\mu$, we say that $E$ is a \emph{leaf} if  the valency of either $\nu$ or $\mu$ is equal to one, 
and we call a vertex of valency one a \emph{leaf vertex}. 

To a Brauer graph $\Gamma$ we associate a quiver $Q_\Gamma$ where the vertices of $Q_\Gamma$ correspond to the edges in $\Gamma$. Let $i$ and $j$ be two distinct vertices in $Q_\Gamma$
 corresponding to edges $E_i$ and $E_j$ in $\Gamma$. Then there is an arrow $i \stackrel{\alpha }{\longrightarrow} j$ in $Q_\Gamma$ if the edge
 $E_j$ is a direct successor of the edge $E_i$ in the cyclic ordering around a common vertex in $\Gamma$. 
If the edge $E_i$ is a leaf with leaf vertex $\nu$ with $e(\nu) \geq 2$ then $E_i$ is its own successor in the cyclic ordering and there is an arrow $i \stackrel{\alpha }{\longrightarrow} i$. 
If $e(\nu) = 1$ then no such arrow exists. 

It follows from the construction of $Q_\Gamma$ that every vertex $\nu$ of $\Gamma$ gives rise to an oriented cycle in $Q_\Gamma$,  
unless $\nu$ is a leaf with leaf vertex $\nu$ where $e(\nu)=1$. 
Furthermore, no two cycles in $Q_\Gamma$ corresponding to distinct vertices in $\Gamma$ have a common arrow. 

Define on $Q_\Gamma$ a set of relations $\rho_\Gamma$ as follows. Let $E_i$ be an edge in $\Gamma $ with vertices $\nu$ and $\nu'$ where  
if $\nu$ (respectively $\nu'$) is a leaf vertex then  
$e(\nu) \neq 1$ (respectively $e(\nu') \neq 1$). Denote by $C_{\nu,i}$ and $C_{\nu',i}$ the corresponding cycles in $Q_\Gamma$ 
starting at vertex $i$ in $Q_\Gamma$. Let $C^{e(\nu)}_{\nu,i}$ be
the $e(\nu)$-th power of $C_{\nu,i}$. Then
 $C^{e(\nu)}_{\nu,i} - C^{e(\nu')}_{\nu',i} \in \rho_\Gamma$. 
Now suppose the edge $E_i$ is such that $\nu$ is a leaf with $e(\nu)=1$. Then $C_{\nu',i} \alpha \in \rho_\Gamma$ where $\alpha$ is the arrow in $C_{\nu',i}$ 
starting at $i$. 
Finally if $\alpha, \beta$ are two arrows in $Q_\Gamma$ such that $t(\alpha) = s(\beta) =i $ where $\alpha$ is in $C_{\nu,j}$ for $s(\alpha) = j$ and 
$\beta$ is in $C_{\nu',i}$ with $\nu \neq \nu'$ then $\alpha \beta \in \rho_\Gamma$. 

The algebra $B_\Gamma = kQ_\Gamma /I_\Gamma$ where $I_\Gamma$ is the ideal generated by $\rho_\Gamma$ is called the \textit{ Brauer graph algebra associated to the Brauer graph} $\Gamma$. 
Note that $B_\Gamma$ is a finite dimensional symmetric special biserial algebra, that is $B_\Gamma$ satisfies (S1) and (S2). Furthermore, it immediately follows from the definition of  Brauer graph algebras that they satisfy condition (S4).

We remark that the above definition of Brauer graph algebras has to be slightly adjusted to include the algebras $k[x]/x^2$ and $k$. The Brauer graph 
of $k[x]/x^2$ would then correspond to a single edge with both vertices of multiplicity one and that of the algebra $k$  to a single vertex.
However, in order to keep the above notation and for clarity of exposition
we exclude these two algebras in what follows. 


 
\subsection{Symmetric special biserial algebras}


Let $\Lambda = kQ/I$ be a finite dimensional symmetric special biserial algebra. 
Since $k$ is algebraically closed, without loss of generality we can assume that a set of relations $\rho$ generating $I$ contains only zero relations and commutativity 
relations of  the form $p-q$ for $p, q$  paths in $Q$ such that $p,q \notin \rho$. 

Since $\Lambda$ is symmetric special biserial the projective indecomposable modules are uniserial or biserial. We adopt the following notation. 
Let $i$ be a vertex in $Q$. Then if the projective 
indecomposable $P_i$ at $i$ is uniserial there exists a unique non-trivial maximal path $p$ in $(Q,I)$ 
with $s(p) = t(p) = i$. 
Let $e_i$ be the trivial path at $i$. 
Then we write $P_i = P_i(p, e_i) = P_i(p)$ or $P(p)$ for short. If $P_i$ is biserial then
there exist two distinct non-trivial paths $p, q$ in $Q$ with $s(p)= s(q) = t(p) = t(q) = i$ such that $p-q \in I$. 
 We write $P_i = P_i(p,q)$ or $P(p,q)$ for short. Since $\Lambda$ is symmetric, the projective indecomposable at $i$ is also the 
injective indecomposable at $i$
 and hence the paths $p$ and $q$ are maximal in $(Q, I)$.  
It is a direct consequence of (S2) that $p$ and $q$ do not start or end with a common arrow. 
Furthermore, for any two non-trivial paths $p,q$ in $Q$, if there is a relation $p-q \in I$ where $p,q \notin I$ then 
$P_i = P(p,q)$ where $i=s(p)$. This follows directly from (S2) and the fact that $\rad P_i / \soc P_i$ is a direct sum of two uniserial modules. 

\begin{remark}
Suppose that $B$ is a Brauer graph algebra with Brauer graph $\Gamma$. With the notation above let $C^{e(\nu)}_{\nu,i}, C^{e(\nu')}_{\nu',i}$  
be the two  maximal paths defined by 
two vertices $\nu, \nu'$ of $\Gamma$ connected by an edge $E_i$ corresponding to a vertex $i \in Q_\Gamma$. Then the projective indecomposable $B$-module 
at $i$ is given by $P_i(C^{e(\nu)}_{\nu,i}, C^{e(\nu')}_{\nu',i})$ where $C^{e(\nu)}_{\nu,i}$ or $ C^{e(\nu')}_{\nu',i}$ might be the trivial path at $i$
 if $\nu$ (respectively $\nu'$) is a leaf vertex with $e(\nu)=1$ (respectively $e(\nu') = 1$). In the latter case the projective indecomposable is uniserial. 
\end{remark}

\begin{Lemma}\label{EqualProjectiveIso}
Let $\Lambda = kQ_\Lambda / I_\Lambda$ and $\Delta = kQ_{\Delta} / I_{\Delta}$ be symmetric special biserial algebras. Suppose $Q=Q_\L = Q_{\D}$ 
and suppose that for every vertex in $Q$, the projective indecomposable modules of $\Lambda$ and $\Delta$ 
have  $k$-bases given by the same paths in $Q$. Then the algebras $\L$ and $\D$ are isomorphic. \end{Lemma}

{\it Proof:} We will show that $I_\L = I_\D$. Denote the projective indecomposable module 
of $\L$ (respectively $\D$) at vertex $i$ by $P_i^\L$ (respectively $P_i^\D$). It is enough to show that
every  generating 
relation for $I_\L$ must also be a generating relation for $I_\D$. Suppose
 $\alpha_1 \ldots \alpha_n$ is a path in $Q$ with $\alpha_1 \ldots \alpha_n \in I_\L$. Now assume that  $\alpha_1 \ldots \alpha_n \notin I_\D$. 
 Then there is a path $p_1$ in $Q$ such that   $p = \alpha_1 \ldots \alpha_n p_1$  and $P_i^\D = P(p,q)$ 
 for some possibly trivial path $q$ and where $i = s(\alpha_1)$. But $p \in I_\L$ since $\alpha_1 \ldots \alpha_n \in I_\L$. 
  Therefore  $P_i^\L \not\simeq P(p,q)$, a contradiction and thus $\alpha_1 \ldots \alpha_n \in I_\D$. By symmetry of the argument
   this implies that $\alpha_1 \ldots \alpha_n \in I_\L$ if and only if $\alpha_1 \ldots \alpha_n \in I_\D$. 
  Now suppose that  $p,q$ are paths in $Q$ with $p,q \notin I_\L$  and $0 \neq p-q \in I_\L$. 
  By (S2) $p$ and $q$ do not start with the same arrow (since otherwise $p=q$). But then the projective $P_i^\L$ with $i = s(p) = s(q)$ is biserial. 
  Hence   $P_i^\L \simeq P(p, q)$. 
  But  then $P(p, q) \simeq P_i^\D$ and thus $p-q \in I_\D$. \hfill $\Box$

 \subsubsection{Graph of a symmetric special biserial algebra}\label{GraphSymmetricSpecialBiserial}
  
In \cite{Ro} Roggenkamp 
showed that if the quiver of a symmetric biserial algebra $\Lambda$ has no double arrows then $\Lambda$ is a Brauer graph algebra. We will adapt Roggenkamp's construction of the Brauer
graph associated to a symmetric special biserial algebra to show that this result holds also for quivers with double arrow, thus holding in general. 

Let $\Lambda = kQ/I$ be a symmetric special biserial algebra. 
We now use the projective indecomposable $\Lambda$-modules to define a graph $G_\Lambda$ with a cyclic ordering of the edges around each vertex. As mentioned
above we follow Roggenkamp's construction for this. However, instead of using the Loewy structure of projective indecomposables as in \cite{Ro}, we consider the arrows and paths
defining the projective indecomposables. This eliminates any ambiguity arising from the existence of double arrows. 
Recall that the projective indecomposable $\Lambda$-modules are either uniserial or biserial and that they are denoted by 
$P_i(p)$ and $P_i(p,q)$ respectively where $p, q$ are maximal cyclic paths from a vertex $i$ in $Q$ to itself. 
For the purpose of this construction we consider the trivial path $e_i$ at vertex $i$
to be a maximal cyclic path at $i$ if the projective indecomposable at $i$ is uniserial.  
This implies that there are 
two maximal cyclic paths at every vertex (one of them possibly being the trivial path). For  non-trivial 
maximal cyclic paths $p$ and $q$ at vertex $i$ there are strictly positive integers $m$ and $n$ such that
 $p = p_0^m$ and $q = q_0^n$ where $p_0$ and $q_0$ are   cycles starting and 
ending at $i$ which are no proper power of cycles of shorter length. We call $p_0$ and $q_0$ \emph{simple cycles}.
 Set $e(p)=m$ and $e(q)=n$ and call it the \emph{multiplicity} of $p$ and respectively of $q$. Trivial maximal paths always have multiplicity one.

Suppose $p  = p_0^m$ is a maximal path as above where $p_0 = i \stackrel{\alpha_0}{\longrightarrow} i_1 \stackrel{\alpha_1}{\longrightarrow} i_2 \cdots \stackrel{\alpha_{k-1}}{\longrightarrow}
 i_k \stackrel{\alpha_{k}}{\longrightarrow}
 i$. Define the \emph{$p$-cycle of $p$} to be the sequence $\mu (p) = (i, i_1, i_2, \ldots, i_k)$ of vertices in $Q$. If the trivial path $e_i$ at a vertex $i$ is maximal, we set
 $\mu (e_i)= i$.

\begin{remark} (1) If $\tilde{p}$ is a cyclic rotation of $p$ then $e(\tilde{p}) = e(p)$. 

(2) The vertices of $Q$ occurring in $\mu(p)$ need not all be different. However, since $\Lambda$ is special biserial, each one can occur no more than twice.

(3) If $\tilde{p}$ is a cyclic rotation of $p$ then $\mu(\tilde{p})$ is a cyclic rotation of $\mu(p)$. 
\end{remark}
  
Set $p \sim \tilde{p}$ if $\tilde{p}$ is a cyclic rotation of $p$. This defines an equivalence relation  on the set of cyclic rotations of $p$. Denote the equivalence class of $p$ by 
$\nu(p)$ and call it a \emph{vertex}. For later use, we also introduce the following closely related terminology, if $p = p_0^{e(p)}$ with $p_0$ a simple cycle then we call the rotation class of $p_0$ a \emph{vertex cycle}. 

The set of vertices $V$ of $G_\Lambda$ consists of the equivalence classes $ V= \{ \nu(p) \; | \;  p \mbox{ maximal path in } (Q,I) \}$. 
\sloppy To each vertex $\nu(p)$ we associate its multiplicity $e(p)$. Note that this 
is well-defined since $e(p) = e(\tilde{p})$ for $p \sim \tilde{p}$. To the vertex $\nu(p)$ we now attach \emph{germs of edges} labelled by the vertices of $Q$ contained
 in $\mu(p)$. Note that the arrows in $p_0$
induce a linear order on the vertices in $\mu(p)$ by setting $i_n < i_{n+1}$  for $ 0 \leq n \leq k-1$ and we can complete this to a cyclic order by setting $i_k < i$. 
This  defines a cyclic
ordering of the germs around $\nu(p)$.  
Call this the \emph{local structure of $G_\Lambda$}.

\begin{remark} (1) The local structure of $G_\Lambda$ is not changed if we replace $\mu(p)$ by $\mu(\tilde{p})$ where $p \sim \tilde{p}$. Hence denote the $p$-cycle $\mu(p)$ by $\mu_\nu$
where $\nu = \nu(p)$. 

(2) In the set $\{ \mu_\nu | \nu \in V\}$ each vertex of $Q$ appears exactly twice. Thus for every vertex in $Q$ there are exactly two germs labelled by that vertex. Note that a vertex of $Q$ can appear twice in one 
$\mu_\nu$, labelling two distinct germs. 
\end{remark}

\begin{Definition}
The graph $G_\Lambda$ of a symmetric special biserial algebra 
$\Lambda$ is given by the vertices $\nu \in V$ with  cyclic ordering given by $\mu(p)$ for some maximal path $p$ such that  $\mu(p) = \mu_\nu$. 
There is an edge from $\nu$ to $\nu'$ if 
the cycles $\mu_\nu$ and $\mu_{\nu'}$ have a common vertex of $Q$. In this case we join the two corresponding germs to a genuine edge. 
\end{Definition}

Note that if a vertex of $Q$ occurs twice in $\mu_\nu$ then there is a loop in $G_\Lambda$. 

The graph $G_\Lambda$ is a graph with a cyclic ordering of the edges at each vertex. Thus together with the multiplicity of the maximal paths (defined above) associated to the
corresponding vertices
of $G_\Lambda$, the graph $G_\Lambda$ is a Brauer graph.
Let $B$ be the corresponding Brauer graph algebra.  
It follows immediately from the definition of the cyclic ordering in $G_\Lambda$ which is induced by the arrows in $Q$ 
that the projective indecomposable $B$-modules have a $k$-basis given by the same paths in $Q$ as the corresponding 
 projective indecomposable $\Lambda$-modules. 
Combining this with Lemma~\ref{EqualProjectiveIso} proves Theorem~\ref{SymmetricSpecialBiserialIsBGA}  which for the convenience of the reader we restate here.

\textbf{Theorem~\ref{SymmetricSpecialBiserialIsBGA}} \textit{
Let $k$ be an algebraically closed field and let $\Lambda = kQ/I$ be a symmetric special biserial algebra and let $G_\Lambda$ be its graph with local structure as defined above. 
Let $B$ be the symmetric Brauer graph algebra associated to $G_\Lambda$, then $B$ is isomorphic to $\Lambda$. 
}



\section{Trivial extensions of gentle algebras}\label{TA}


Let $A$ be a finite dimensional  $k$-algebra.  Recall that throughout we assume that the field $k$ is algebraically closed. The\textit{ trivial extension} $T(A)$ of $A$ by its minimal injective co-generator $D(A)$ is  the algebra
$T(A)=A \ltimes D(A)$. As a $k$-vector space $T(A)$ is given by $A \oplus D(A)$  with the multiplication
defined by $(a, f)(b, g) = (ab, ag + fb)$, for $a, b \in A$ and $f,g \in D(A)$. It is well-known that the trivial extension algebra $T(A)$
is  a symmetric algebra (see, for example, \cite[Proposition 6.5]{Schi}).

As recalled in the introduction  it is shown in \cite{PS, Ringel} and \cite{Schroer} that 
  $A$ is gentle if and only if  $T(A)$ is special biserial. 
Therefore the trivial extension of a gentle algebra is a  symmetric special biserial algebra. 
Hence it follows from Theorem~\ref{SymmetricSpecialBiserialIsBGA}  
that $T(A)$ is a Brauer graph algebra.

Let $A = kQ_A/I_A$. It is proven in \cite[2.2]{FP1} that the vertices of the quiver $Q_{T(A)}$ of $T(A)$ correspond to 
the vertices of the quiver $Q_A$ of $A$ and 
that the number of arrows from a vertex $i$ to a vertex $j$  in $T(A)$ is equal to the number of arrows from $i$ to $j$ in $Q_A$ plus 
the dimension of the $k$-vector space $e_j ({\rm soc}_{A^e} A) e_i$. 
 
The map $D(A) \longrightarrow T(A)$ given by $x \mapsto (0,x)$ is an injective $T(A)$-bimodule homomorphism resulting in a short exact sequence of $T(A)$-bimodules
$$ 0 \longrightarrow D(A) \longrightarrow T(A) \longrightarrow A \longrightarrow 0.$$
For every vertex $i$ in $Q_{T(A)}$, let $e_i \in T(A)$ be the corresponding primitive idempotent. Then by left multiplication with $e_i$ we obtain a short exact sequence of right $T(A)$-modules
$$ 0 \longrightarrow e_i D(A) \longrightarrow e_i T(A) \longrightarrow e_i A \longrightarrow 0$$ where $e_iT(A)$ is the projective-injective indecomposable $T(A)$-module at vertex $i$, $e_iA$ is the projective $A$-module at $i$ and $e_i D(A)$ is the injective $A$-module at $i$ (see, for example, \cite[Section 6.2]{Schi}). 

\subsection{Construction of the graph of a gentle algebra}

Given an indecomposable gentle algebra $A = kQ/I$, we will define a graph $\Gamma_A$ associated to $A$. 

As in the construction of the graph associated to a symmetric special biserial algebra, we will now extend the set $\M$ of maximal paths in $(Q,I)$. Let $\oM$ be the set 
containing $\M$ and all those trivial paths associated to  vertices $i$ of $Q$ that are a sink with  a single arrow ending at $i$,  a source with   
 a single arrow starting from $i$ or such that there is  a single arrow $\alpha$ ending at $i$ and a single arrow $\beta$ starting at $i$ and  $\alpha \beta \notin I$. 
 We will still call the elements of $\oM$ maximal paths.  
 As a consequence of the above definition every vertex of $Q$ lies in exactly two distinct maximal paths in $\oM$.

Let $\Gamma_A$ be the graph with vertex set $V_0$ given by the (extended) set of maximal paths $\oM$, and with edge set
$V_1$ given by the set of vertices of $Q$. Denote by $\nu(m)$ the vertex of $\Gamma_A$ corresponding 
to the maximal path $m$ in $\oM$ and denote by $E_i$ the edge in $\Gamma_A$ corresponding to the vertex $i$ in $Q$. 
Recall that by the definition of $\oM$ every vertex of $Q$ lies in exactly two maximal paths. 
Then for a vertex $i$ in $Q$ lying in the maximal paths $m_i$ and $n_i$ of $\oM$, let $E_i$ be the edge in $\Gamma_A$ connecting 
$\nu(m_i)$ and $\nu(n_i)$.  




\begin{Lemma}\label{LinearOrder} Each  maximal path $m$ in $\M$ gives rise to a linear order of the edges connected to the corresponding vertex $\nu(m)$ in $\Gamma_A$. \end{Lemma} 

{\it Proof:}  Let $m = i_1 \stackrel{\alpha_1}{\longrightarrow} i_2 \stackrel{\alpha_2}{\longrightarrow} \cdots
 \stackrel{\alpha_n}{\longrightarrow} i_{n+1}$ be a maximal path in $\M$ and let $\nu(m)$ be the corresponding 
 vertex in $\Gamma_A$. Then each $i_j$ corresponds to an edge $E_{i_j}$ in $\Gamma_A$ connecting $\nu(m)$ and 
 $\nu(n_j)$, where 
 $n_j$ is the second maximal path in $\oM$ going through the vertex $i_j$ of $Q$. 
Thus as edges in $\Gamma_A$, all $E_{i_j}$, for $1 \leq j \leq n+1$, are connected to $\nu(m)$. 
Then setting $E_{i_j} < E_{i_l}$ if there exists a subpath $p$ of $m$ with $s(p) = i_j$ and $t(p) = i_l$  
defines a linear order on the set $\{E_{i_1}, E_{i_2}, \ldots, E_{i_{n+1}}\}$.
  \hfill $\Box$

Note that maximal paths in $A$ correspond to maximal fans in $\Gamma_A$. A fan is a subgraph of $\Gamma_A$ such that all edges in the subgraph are connected to a common vertex. 

An example of a graph associated to a gentle algebra is given by the arcs of a triangulation of a marked Riemann surface: Let $T$ be an
ideal triangulation of a bordered
 Riemann surface with marked points in the boundary and let $A$ be the associated gentle algebra \cite{ABCP}. Then $\Gamma_A = T$ (see Section~\ref{AdmissibleCutsSection} for more details and an example).

We will now construct a quiver $Q_E$ extending the quiver $Q_A$. It follows from the proof of Lemma~\ref{LinearOrder}  that if $i_1 < i_2 < \cdots < i_{n+1}$ is the linear order defined by 
a maximal path $m \in \M$ then this linear order can be completed to a cyclic order by adding a single arrow  $i_{n+1} \stackrel{\beta_m}{\longrightarrow} i_1$. 

Define $Q_E$ to be the quiver with 
\begin{itemize}
\item set of vertices  given by the vertices in $Q$ and 
\item set of arrows given by the arrows of $Q$ together with a set of new arrows 
$\{\beta_m | \mbox{ for every $m \in \M$ where $s(\beta_m) = t(m)$ and $t(\beta_m)= s(m)$} \}$.
\end{itemize}   

\begin{Lemma}
With the notations above the quiver $Q_E$ is the quiver of the trivial extension algebra $T(A)$ of $A$. 
\end{Lemma}

{\it Proof:}   By \cite{FP1} and Lemma~\ref{MaximalPathsBasis}  the arrows of $Q_{T(A)}$ are given by the arrows of $Q_A$ plus for every maximal path $m \in \M$   a new arrow $\beta_m$ 
 such that $s(\beta_m) = t(m)$ and $t(\beta_m)= s(m)$.  \hfill $\Box$

The following Lemma follows directly from the definitions of $Q_E$ and the definition of the quiver of a Brauer graph algebra. 

\begin{Lemma}\label{SameQuiver}
With the notations above the quiver $Q_E$ is the quiver of the Brauer graph algebra with associated Brauer graph $\Gamma_A$ and with cyclic ordering of the edges of $\Gamma_A$ around each vertex  induced by the arrows in $Q_E$. 
\end{Lemma}

We can now prove Theorem~\ref{TrivialExtensionofGentleAlgebra} which we restate here for convenience. 

\textbf{Theorem~\ref{TrivialExtensionofGentleAlgebra}} \textit{
Let  $A= kQ/I$ be a gentle algebra over an algebraically closed field $k$ and let $\Gamma_A$ be its graph. Let $B$ be the Brauer graph algebra defined on $\Gamma_A$ with cyclic ordering induced by the maximal paths in $A$ and 
with multiplicity equal to one at every vertex. Then the trivial extension $T(A)$ of $A$ is isomorphic to $B$ and  $G_{T(A)} = \Gamma_A$ where $G_{T(A)}$ is the graph of $T(A)$ as symmetric special biserial algebra. }

{\it Proof:} The Brauer graph algebra associated to  $\Gamma_A$  and the algebra $T(A)$ are  symmetric special biserial algebras. The projective indecomposable modules  of the Brauer graph algebra can be read of the graph $\Gamma_A$ as described above.  The projective indecomposable $T(A)$-modules are given by short exact sequences as described in Section~\ref{TA} above. That is   they are induced by the short exact sequences
$$ 0 \longrightarrow e_i D(A) \longrightarrow e_i T(A) \longrightarrow e_i A \longrightarrow 0$$ 
 for every vertex $i$ in $Q$. Suppose that $i$ is a vertex such that there are two non-trivial maximal paths 
 $m_i$ and $n_i$ of $A$ going through $i$ where
 $$m_i = i_0 \stackrel{\alpha_0}{\longrightarrow} i_1 \stackrel{\alpha_1}{\longrightarrow}  \cdots \stackrel{\alpha_{k-1}}{\longrightarrow} i_k = i \stackrel{\alpha_k}{\longrightarrow} \cdots
 \stackrel{\alpha_{m-1}}{\longrightarrow} i_m$$
and 
$$  n_i = j_0 \stackrel{\gamma_0}{\longrightarrow} j_1 \stackrel{\gamma_1}{\longrightarrow}  \cdots \stackrel{\gamma_{l-1}}{\longrightarrow} j_l = i \stackrel{\gamma_l}{\longrightarrow} \cdots
 \stackrel{\gamma_{n-1}}{\longrightarrow} j_n.$$
 
 The projective right $A$-module $e_i A$ has a $k$-basis given by all paths starting at $i$. By (S1) there are at most two arrows starting at $i$, which are $\alpha_k$ and $\gamma_l$ (if they exist). By (S2) there is at most one arrow $\beta$ such that $\alpha_k \beta$ is non-zero in $A$, so $\beta = \alpha_{k+1}$ and there is at most one arrow $\delta$ such that $\alpha_{k+1} \delta$ is non-zero, so $\delta = \alpha_{k+2}$ and so forth, until we reach $\alpha_{m-1}$. Then by maximality of $m_i$ there exists no arrow $\varepsilon$ in $Q$ such that $\alpha_{m-1} \varepsilon$ is non-zero in $A$. Similarly, $ i \stackrel{\gamma_l}{\longrightarrow} \cdots
 \stackrel{\gamma_{n-1}}{\longrightarrow} j_n$ is the unique non-zero path starting with $\gamma_l$ and there exists no arrow $\beta$ in $Q$ such that $\gamma_l \ldots \gamma_{n-1} \beta$ is non-zero in $A$. Therefore, $e_iA$ corresponds to the string module $M(w)$ with string $w = i_m \stackrel{\alpha_{m-1}}{\longleftarrow}  \cdots  \stackrel{\alpha_{k+1}}{\longleftarrow} i_{k+1} \stackrel{\alpha_{k}}{\longleftarrow} i  \stackrel{\gamma_l}{\longrightarrow} j_{l+1} 
 \stackrel{\gamma_{l+1}}{\longrightarrow} \cdots \stackrel{\gamma_{n-1}}{\longrightarrow} j_n$. Note that if neither the arrow $\alpha_k$ nor the arrow $\gamma_l$ exists then $e_iA$ is the simple module at $i$. 
  Similar arguments on the arrows ending at $i$ show that  $e_i D(A)$  is isomorphic to the string module $M(v)$ with string $v= i_0 \stackrel{\alpha_0}{\longrightarrow} i_{1} \stackrel{\alpha_1}{\longrightarrow}  \cdots  \stackrel{\alpha_{k-2}}{\longrightarrow}
 i_{k-1}\stackrel{\alpha_{k-1}}{\longrightarrow} i \stackrel{\gamma_{l-1}}{\longleftarrow}    j_{l-1} \stackrel{\gamma_{l-2}}{\longleftarrow} \cdots  \stackrel{\gamma_{1}}{\longleftarrow} j_1 \stackrel{\gamma_{0}}{\longleftarrow} j_0$. In $Q_{T(A)}$ there are  
 arrows $i_m  \stackrel{\beta_{m_i}}{\longrightarrow} i_0$ and  $j_n  \stackrel{\beta_{n_i}}{\longrightarrow} j_0$ such that the projective-injective indecomposable $T(A)$-module
  $e_i T(A)$ is given by 
 the biserial module that has the simple at $i$ as top and socle and whose heart $\rad e_iT(A) / \soc  e_i T(A) $ is given by a direct sum of uniserial modules given by the direct strings
 $i_{k+1} \stackrel{\alpha_{k+1}}{\longrightarrow}   \cdots  \stackrel{\alpha_{m-1}}{\longrightarrow} i_m \stackrel{\beta_{m_i}}{\longrightarrow} i_0
  \stackrel{\alpha_0}{\longrightarrow} i_1 \stackrel{\alpha_1}{\longrightarrow}  \cdots \stackrel{\alpha_{k-2}}{\longrightarrow} i_{k-1}$ and 
  $j_{l+1} \stackrel{\gamma_{l+1}}{\longrightarrow}   \cdots  \stackrel{\gamma_{n-1}}{\longrightarrow} j_n \stackrel{\beta_{n_i}}{\longrightarrow} 
  j_0 \stackrel{\gamma_0}{\longrightarrow} j_1 \stackrel{\beta_1}{\longrightarrow}  \cdots \stackrel{\beta_{l-2}}{\longrightarrow} j_{l-1}$.
   Denote by $E_i$ the edge of $\Gamma_A$ corresponding to the vertex $i$ of 
  $Q$ and let $\nu(m_i)$
   and $\nu(n_i)$ be the vertices at either end of $E_i$ corresponding to the maximal paths $m_i$ and $n_i$. We have seen that the arrows in $m_i$ and $n_i$
    give rise  to a cyclic ordering of the edges $i_s$ around $\nu(m_i)$ and $j_t$ around $\nu(n_i)$. Then the 
   projective-injective indecomposable $B$-module 
  $P^B_i$ at vertex $i$  has a  $k$-basis given by   paths in $Q_E$ and $e_i T(A)$ has a $k$-basis given by the same paths in $Q_E$.  
  
  We will now consider the situation  where either $n_i$ or $m_i$ correspond to the trivial path at $i$. Assume without loss of generality that $n_i = e_i$. Then $e_i T(A) $ is uniserial corresponding to the string 
  $ i_k = i   \stackrel{\alpha_k}{\longrightarrow} i_{k+1} \stackrel{\alpha_{k+1}}{\longrightarrow}   \cdots  \stackrel{\alpha_{m-1}}{\longrightarrow} i_m 
  \stackrel{\beta_{m_i}}{\longrightarrow} i_0 \stackrel{\alpha_0}{\longrightarrow} i_1 \stackrel{\alpha_1}{\longrightarrow}  \cdots \stackrel{\alpha_{k-1}}{\longrightarrow} i_{k} = i$ and the  edge $E_i$ in $\Gamma_A$ corresponding to $i$ is a leaf.
   By an argument similar to the one above we see that $P_i^B$, the projective indecomposable $B$-module at 
   vertex $i$,  and $e_i T(A)$ both have a basis given by the same paths in $Q_E$. 
   
   By Lemma~\ref{EqualProjectiveIso}, the algebras $T(A)$ and $B$ are isomorphic. 
  From the structure of the projective indecomposable modules for $B$ and $T(A)$ it immediately follows that the graphs $\Gamma_A$ and $G_{T(A)}$ are isomorphic. That the 
  cyclic orderings of the edges around each vertex in  $\Gamma_A$ and $G_{T(A)}$ coincide, follows from the fact that both cyclic orderings are induced by the arrows in $Q_E$. 
    \hfill $\Box$






\section{ Admissible Cuts in Symmetric Special Biserial Algebras}\label{AdmissibleCutsSection}

One way of deleting arrows in quivers and constructing new algebras in this way is through the notion of admissible cuts of finite dimensional algebras which  has  been studied 
 for example in \cite{BFPPT, F, FP2},  and in the context of cluster algebras, see for example, \cite{DS},  and  also \cite{MP}.

Given a Brauer graph algebra $\Lambda = kQ_\Lambda/I_\Lambda$ with multiplicity one at all vertices of the corresponding Brauer graph, we define a \emph{  cutting set $D$ of $Q_\Lambda$}
 to be a subset of arrows
in $Q_\Lambda$ formed of exactly one arrow in every vertex cycle in $Q_\Lambda$ 
  (see section~\ref{GraphSymmetricSpecialBiserial} for the definition of a  vertex cycle). Note that this corresponds  to the cutting set defined in \cite{FP2}.
An \emph{ admissible cut of $\Lambda$ with respect to the cutting set $D$} 
is the algebra $kQ_\Lambda / J_\Lambda$ where $J_\Lambda$ is the ideal generated by $I_\Lambda  \cup D$.

 
We now prove Theorem~\ref{ArrowDeletionTheorem} which  is in a similar vein as the results in \cite{FP2} (see also \cite{F}) and which states: 
 
 \textbf{Theorem~\ref{ArrowDeletionTheorem}} \textit{
Let $k$ be an algebraically closed field and let $\Lambda = kQ_\Lambda / I_\Lambda$ be a  Brauer graph algebra with multiplicity one at all vertices in the associated Brauer graph.
 Let $A$ be an admissible cut of $\Lambda$. 
Then $A$ is gentle and $T(A)$ and $\Lambda$ are isomorphic. 
}

{\it Proof: }  Let $Q$ be the quiver such that $Q_0 = (Q_\Lambda)_0$ and $Q_1 = (Q_\Lambda)_1 \setminus D$ where $D$ is a cutting set of $Q_\Lambda$. 
Define $I $ to be the ideal $I_\Lambda \cap kQ$ and let $J_\Lambda$ be the ideal of $\Lambda$ generated by $I_\Lambda \cup  D$.
 Then it follows from the second isomorphism theorem that $A= kQ_\Lambda / J_\Lambda \simeq kQ /I$.

We first show that $kQ/I$ is gentle. Since $Q \subset Q_\Lambda$, it clearly satisfies (S1). 
All zero relations in $I_\Lambda$ are of the form $p-q$ and $p\alpha$ for $p,q$ simple cycles and  $\alpha$ an arrow in $Q_\Lambda$ or are 
 monomial relations $\alpha \beta$ of lengths two where $\alpha$ and $\beta$ belong to two distinct vertex cycles. 
By definition of $D$, any maximal cyclic path $p$  appearing in a relation generating $I_\Lambda$, contains exactly one arrow in $D$. 
Thus it is not in $I = I_\Lambda \cap kQ$ and therefore $I$ is generated by paths of length two only. 
In order to show that (S4) holds, suppose that $\alpha, \beta, \beta'$ are  arrows in $Q \subset Q_\Lambda$ such that $\alpha \beta \in I$ and $\alpha \beta' \in I$. Then 
$\alpha \beta \in I_\Lambda \supset I$ and also $\alpha \beta' \in I_\Lambda \supset I$ and $\beta = \beta'$, since (S4) holds for $\Lambda$.  
A similar argument holds for arrows preceding $\alpha$ and thus (S4) holds for $kQ/I$. 
To show (S2),  suppose $\alpha, \beta, \beta'$ are arrows in $Q$ such that $\alpha \beta$ and $ \alpha \beta'$ are non-zero in $kQ$. Suppose further that  $ \alpha \beta \notin I$
and $ \alpha \beta' \notin I$.  
Since $I = I_\Lambda \cap kQ$ this implies $ \alpha \beta  \notin I_\Lambda $ and $ \alpha \beta' \notin I_\Lambda $.
Because (S2) holds for $\Lambda$ this 
implies $\beta = \beta'$. By a symmetric argument for arrows preceding $\alpha$, it follows that (S2) holds for $kQ/I$ and therefore   $A \simeq kQ/I$ is gentle.   

There is a bijection between vertex cycles in $\Lambda$ and maximal paths in $kQ/I$. Namely let $i \stackrel{\alpha_0}{\longrightarrow} i_1 \stackrel{\alpha_1}{\longrightarrow} \ldots \stackrel{\alpha_{n-1}}{\longrightarrow} i_n \stackrel{\alpha_n}{\longrightarrow} s$
be a vertex cycle $\nu$ in $\Lambda$. Suppose $\alpha_j \in D$ for some $0 \leq j \leq n$. Then no other arrow in $\nu$ is in $D$ and
 $ p = i_{j+1}  \stackrel{\alpha_{j+1}}{\longrightarrow} i_{j+2} \ldots 
i_n \stackrel{\alpha_n}{\longrightarrow} i  \stackrel{\alpha_0}{\longrightarrow} i_1 \ldots i_{j-1} \stackrel{\alpha_{j-1}}{\longrightarrow} i_j$ is a path in $(Q,I)$.  As 
subpaths of $p$, we have $\alpha_{j-1}
\alpha_j \notin I$ and $\alpha_j \alpha_{j+1}  \notin I$. If there exists an arrow $\beta$ in $Q$ with $s(\beta) = t(\alpha_{j-1})$ then by (S2)  $\alpha_{j-1} \beta \in I$. And hence $p\beta \in I$. Similarly, if there
exists an arrow $\gamma \neq \alpha_{j}$ in $Q$ such that $t(\gamma) = s(\alpha_{j+1})$ then $\gamma p \in I$.   Hence 
$p$ is a maximal path in $(Q,I)$.  Conversely,  every non-zero path in $\Lambda$ is a  subpath of a unique vertex cycle. 
Thus if $p$ is a maximal path in $kQ/I$ then $p$ is  a subpath of a vertex cycle $\nu$. But since 
no two distinct maximal paths in $kQ/I$ can  have a common arrow and since we have cut exactly one arrow in each vertex cycle of $\Lambda$, $p$ is the path that starts at the end
of the cut arrow in $\nu$ and ends at the start of the cut arrow in $\nu$. Thus the vertices of the graph $G_\Lambda$ of the symmetric special biserial algebra $\Lambda$ and
the vertices of the graph $\Gamma_{kQ/I}$ of the gentle algebra $kQ/I$ coincide.

Furthermore, the vertices of $Q_\Lambda$ in a $p$-cycle associated to a vertex cycle $\nu$ correspond up to rotation to the vertices of the corresponding 
maximal path $p$ in $(Q,I)$. Recall further that the edges in $G_\Lambda$ are given by connecting the two $p$-cycles 
containing the same vertex of $Q_\Lambda$ and 
that the edges in $\Gamma_{kQ/I}$ are given by connecting the two maximal 
paths containing the same vertex of $Q$. The cyclic order of the edges around a vertex in $G_\Lambda$ is induced by the arrows in the corresponding vertex cycle and the maximal paths in $(Q,I)$ 
induce the linear order of the
edges in $\Gamma_{kQ/I}$. As described above the latter can be extended to a cyclic ordering.  
With this induced cyclic ordering on $\Gamma_{kQ/I}$,  the graphs $G_\Lambda$ and $\Gamma_{kQ/I}$ are isomorphic   and have the same cyclic ordering of the edges around each vertex. Thus by Theorem~\ref{SymmetricSpecialBiserialIsBGA} and Theorem~\ref{TrivialExtensionofGentleAlgebra},
the algebras $T(A) \simeq T(kQ/I) $ and $ \Lambda$ are isomorphic. \hfill $\Box$

Note that a Brauer graph algebra with multiplicity function identically equal to one corresponds to a symmetric special biserial algebra 
$kQ/I$ where all relations in $I$ that are not monomial of length 2, are given by simple cycles.

The following Corollary immediately follows from Theorem~\ref{ArrowDeletionTheorem}.

\textbf{Corollary~\ref{BGAisTrivialExtension}} \textit{
Every Brauer graph algebra over an algebraically closed field with multiplicity function identically one is the trivial extension of a gentle algebra.
}

Finally we end the paper with two applications to gentle algebras associated to marked Riemann surfaces. 

Let $S$ be a connected 2-dimensional Riemann surface with boundary $\d S$ and let $M$ be a non-empty finite set of points in the boundary  of $S$. 
Let $\T$ be an ideal triangulation of $(S,M)$ as defined in \cite{FST}. That is $\T$ consists of a maximal collection of arcs 
given by isotopy classes of pairwise non-intersecting 
curves connecting two marked points and such that the curves are not isotopic to a curve lying in the boundary connecting 
two adjacent marked points on the same boundary component of $S$.  
In \cite{ABCP} and in \cite{Lab}, a finite dimensional algebra, the so-called Jacobian algebra, is associated to the triple $(S,M,\T)$. In \cite{ABCP} this algebra is shown to be gentle.

\textbf{Corollary~\ref{Jacobian}} \textit{
Let $k$ be an algebraically closed field and let the $k$-algebra $A$ be the Jacobian algebra associated to
 an ideal triangulation $\T$ of a marked Riemann surface $(S,M)$ with $M \subset \d S$.
 Then the trivial extension $T(A)$ of $A$ is isomorphic to 
the Brauer graph algebra with Brauer graph $\T$ where the cyclic ordering of the edges around each vertex in $\T$ is induced by the orientation of $S$. 
 }

{\it Proof:}  Let $\Lambda = kQ_\Lambda / I_\Lambda$ be the Brauer graph algebra associated to $\T$ and let $A = kQ_A / I_A$ 
be the (gentle) Jacobian algebra associated to the 
triangulation $\T$ of $(S,M)$. We start by constructing an admissible cut  $J =kQ_J /I_J$ of $\Lambda$. For this embed the quiver $Q_\Lambda$
 into $(S,M,\T)$ via $\T$. That is 
every arc of $T$ corresponds to a vertex of $Q_\Lambda$ and the arrows of $Q_\Lambda$ correspond to the ordering of the edges of 
$\T$ in the cyclic ordering induced by 
the orientation of $S$. 
Then  each vertex of $\T$ that is not a leaf vertex, corresponds to a cycle in $Q_\Lambda$ and there is exactly one arrow in every such cycle
 crossing either one or two boundary segments of $\T$ (crossing one boundary segment
precisely when the  vertex under consideration is the only marked point in its boundary component of $(S,M)$). 
The collection of these arrows is a cutting set $D$ of
$Q_\Lambda$ and $Q_A = Q_\Lambda \setminus D =Q_J$. For an example, see figure~1 below. 

Let $\rho_\Lambda$ be the set of relations generating $I_\Lambda$ as described in Section~\ref{BGA}. Then $\rho_\Lambda \cap kQ_J$ is a generating set of 
relations for $I_J$. Let
 $\rho_A$ be a set of relations generating $I_A$. Then $\rho_A$ consists of paths $\alpha \beta$ of lengths two where $\alpha$ and $\beta$ are two 
 consecutive arrows in an internal
  triangle of $\T$.  As arrows of  $Q_\Lambda$,  $\alpha$ and $\beta$ are in two distinct vertex cycles of $\Lambda$ and $\alpha \beta \in \rho_\Lambda$.  
  Hence, since $Q_A = Q_J$,  we have $\alpha \beta \in  \rho_\Lambda \cap kQ_J $. 
Conversely, any non-zero path $\alpha \beta$ in $\rho_J$  is given by two arrows $\alpha$ and $\beta$ belonging to two distinct vertex cycles of $\Lambda$. 
Thus $\alpha \beta$ 
corresponds to two consecutive arrows in an internal triangle 
of $T$. From this we conclude that the generating sets of $I_J$ and $I_A$ coincide and hence $A$ is isomorphic to $J$. The fact that $T(A)$ is isomorphic to 
$\Lambda$ then follows 
from Theorem~\ref{ArrowDeletionTheorem}.  \hfill $\Box$

\bigskip

\textbf{ Example.} We consider the gentle Jacobian algebra given by the triangulation of the marked annulus in figure~1. 


\begin{center}
\begin{tikzpicture}[scale=2.5, trans/.style={thick,<-,shorten >=2pt,shorten <=2pt,>=stealth}] 
\draw (0,0) circle (1cm) ;
\draw[fill=gray] (0,0) circle  (.3cm) ;
\node (a) at (90:1cm) [scale=.4] {};
\node (b)at (270:1cm) [scale=.4] {};
\node (c) at (90:.3cm) [scale=.4] {};
\path (a.center) edge node [pos=.5, fill=white,outer sep=1mm,scale=.5]{1} (c.center);
\draw (a.center) .. controls +(0,-.1) and +(-1,0) .. ($(a.center)!.8!(b.center)$) .. controls +(.5,0) and +(.5,.5) .. (c.center);
\node  at (-.55,-.2) [fill=white,outer sep=4mm,scale=.5]{2};
\draw (c.center) .. controls +(.4,.8) and +(1,.5) .. (b.center);
\node  at (.55,-.2) [fill=white,outer sep=4mm,scale=.5]{3};

\node (aa) at (90:.7cm) [scale=.4] {};
\node (bb) at (75:.5cm) [scale=.4] {};
\node (cc) at (60:.46cm) [scale=.4] {};
\node (dd) at (102:.717cm) [scale=.4] {};
\draw[trans, color=red] (dd.center) -- (aa.center);
\draw[trans, color=blue, dashed] (aa.center) .. controls (.5,1.4) and (-.5, 1.4) .. (dd.center);
\node[scale=.7, color=blue] at (0,1.3){$\alpha$};
\draw[trans, color=red] (bb.center) -- (aa.center);
\draw[trans, color=red] (cc.center) -- (bb.center);
\draw[trans, color=blue, dashed] (aa.center) .. controls +(-.5,-.5) and +(.1,-.5) .. (cc.center);
\node[scale=.7, color=blue] at (-.245,.35){$\beta$};
\end{tikzpicture}

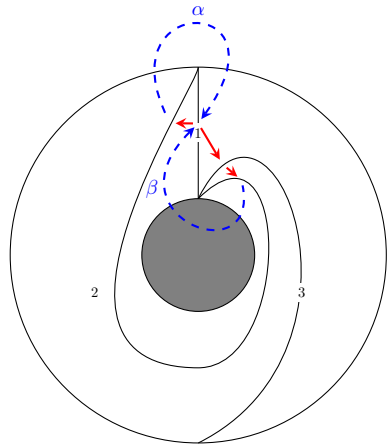
Figure 1: Triangulation of an annulus and associated arrows as described in Corollary 1.5 and its proof. 
\end{center}

The quiver given by the solid (red) arrows is the quiver of the  gentle Jacobian algebra associated to the triangulation $T$ of the annulus in figure 1.  
Note that the vertices of this quiver correspond to the three internal arcs of $\T$, marked 1,2, and 3. The graph given by these arcs together with 
the marked points  constitutes precisely the graph of $\T$. 
The quiver given by the solid (red) arrows together with the dashed (blue) arrows is the quiver of the trivial
 extension of $A$ or equivalently the quiver of the Brauer graph algebra associated to $\T$ as a Brauer graph (with clockwise ordering of the edges
 around each vertex induced by the clockwise orientation of the annulus).  
 The dashed (blue) arrows form the cutting set $D$ described in the proof of Corollary 1.5 above. Note that $\alpha$ corresponds to an arrow crossing two 
 distinct boundary segments of $(S,M)$ and $\beta$ is an arrow crossing a single boundary segment of $(S,M)$.

\bigskip

Given a triangulation $\T$ of   $(S,M)$, in \cite{DS}  surface algebras were defined by cutting $T$ at internal triangles. This gives
rise to a partial triangulation of $(S,M)$. A surface algebra $A$
is again a finite dimensional gentle algebra and it can be constructed by associating a bound quiver $(Q,I)$ to the partial triangulation of $(S,M)$ such 
that $A = kQ/I$ (see \cite{DS}). 

\textbf{Corollary~\ref{SurfaceAlgebra}} \textit{
Let $k$ be an algebraically closed field and let the $k$-algebra $A$ be a surface algebra of a partial triangulation $\P$ of a  Riemann surface $(S,M)$ with set of marked points $M \subset \d S$.
 Then the trivial extension $T(A)$ of $A$ is isomorphic to 
the Brauer graph algebra with Brauer graph $\P$ with cyclic orientation induced by the orientation of $S$. 
}

{\it Proof:} Again, as in the proof of Corollary~\ref{Jacobian},
 we construct  an admissible cut $J$ of the Brauer graph algebra associated to $\P$ and 
we show that $J$ is isomorphic to $A$. However, in $A$ not all relations result from paths of lengths two lying in an internal triangle. There are also relations generated by paths of lengths two lying in a 
quadrilateral resulting from the 'cut' of an internal triangle (see \cite{DS}). However, these relations also appear in the corresponding Brauer graph algebra and its admissible cut. 
Taking this into account, the proof  
then is similar to the  proof of Corollary~\ref{Jacobian}.  \hfill $\Box$


\bibliographystyle{ACM}

\end{document}